\input amstex
\documentstyle{amsppt}

\loadeufb
\loadeusb
\loadeufm
\loadeurb
\loadeusm

\magnification =\magstep 1
\refstyle{1}
\NoRunningHeads

\topmatter
\title On uniformization of compact Kahler manifolds 
\endtitle

\author Robert  Treger \endauthor
\address Princeton, NJ 08540  \endaddress
\email roberttreger117{\@}gmail.com \endemail
\endtopmatter

\document

The aim of the present note is to extend the uniformization theorem of projective manifolds in \cite{9, Introduction, Theorem} to compact Kahler manifolds. In an email  to the author (January, 2015), Dennis Sullivan essentially raised the question whether one can generalize the uniformization theorem in \cite{9}. 
 The author would like to thank him for the question.

Let $X$ be a compact complex manifold of dimension $n\geq 2$. We denote its universal covering by $U_X$. We will derive the following theorem from a similar theorem in \cite{9}.

\proclaim{Theorem (uniformization)} Let $X$ be a compact Kahler manifold of dimension $n$
with  large and residually finite fundamental group $\pi_1(X)$.  If $\pi_1(X)$ is, in addition, nonamenable then $U_X$ is a bounded domain in $\bold C^n$\!. Thus $X$ is projective  by  Poincar\'e  \cite{5,\,Theorem 5.22}.
\endproclaim

\demo{Proof of Theorem} By a  theorem of Moishezon \cite{7}, it will suffice to establish that $X$ is a Moishezon manifold. 
In \cite{4, Sect.\;3},  Gromov uses his notion of Kahler
hyperbolicity to obtain holomorphic $L_2$ forms on $U_X$ and prove that $X$ is 
Moishezon. A priori, we do not know if there are holomorphic $L_2$ forms on $U_X$.

Set $\Gamma: =\pi_1(X)$. Let $\Cal L$ be an arbitrary complex line bundle on $U_X$. We will consider a section $f\in H^0(\Cal L^q,U_X)$ which is not assumed
to be  $L_p$, where $p < \infty$.
As in Koll\'ar \cite{5, Chap\;13.1}, we will employ $\ell^p$ sections $f$  {\it on orbits} of $\Gamma$ in place of $L_p$ sections.
Of course, we need a natural $\Gamma$-invariant Hermitian {\it quasi}-metric on $\Cal 
L^q$ (see the definition in the proof of Lemma 3).

Given an arbitrary $\Gamma$-invariant Hermitian metric on $U_X$, we get the
induced Riemannian metric on $U_X$ with the volume form $d\mu$. Since $\Gamma$ is nonamenable, we get non-constant  bounded harmonic functions 
on $U_X$ by Lyons and Sullivan \cite{6}. 
Employing their theorem, Toledo \cite{8} has established that the space of bounded harmonic functions as well as the space generated by bounded positive harmonic functions are infinite dimensional (see \cite{9, Sect.\;2.6}). 
Given $r$ linearly  independent functions $g_1,\dots,g_r$ on
$U_X$,  clearly  there exist $r$ points $ u_1,\dots,u_r
 \in U_X$ such that the vectors  
$
\langle g_1(u_i),\dots,g_r(u_i)\rangle\; (1\leq i\leq r)
$ 
are linearly independent. 

Let $Har(U_X)$ ($Har^b(U_X)$) be the space of  harmonic functions  (bounded  harmonic functions, respectively) on $U_X$.

In place of the standard $L_2(d\mu)$ integration with the standard Riemannian measure $d\mu$ on $U_X$, we will integrate  the {\it bounded} harmonic functions  with respect to the measure 
$$
dv:=p_{U_X}(s,x,\bold Q) d\mu,
$$
where $\bold Q \in U_X$ is a fixed point  and $p_{U_X}(s,x,\bold Q)$ is the heat kernel. Because all  bounded harmonic  functions are square integrable, i.e. in $L_2(dv)$, we obtain the pre-Hilbert space of bounded harmonic  functions (compare \cite{9, Sect.\;2.4 and Sect.\;4}).
 We observe that the latter pre-Hilbert  space has a completion in the (real)
Hilbert space $H$ of all harmonic  $L_2(dv)$ functions:
$$
H := \biggl\{h \in Har(U_X)\;\biggl |\; \parallel h\parallel^2_H:=  \int_{U_X}
|h (x)|^2 dv < \infty   \biggl\}. 
$$
Let $H^b \subseteq H$ be the Hilbert subspace generated by $Har^b(U_X)$.
These Hilbert spaces are separable infinite dimensional and have  reproducing kernels. 
The group $\Gamma$ acts  isometrically on $H^b: \psi \mapsto
(\psi\circ\gamma)\, (\gamma\in \Gamma)$.

Let $\{\phi_j\}\subset Har^b(U_X)$ be an orthonormal basis of $H^b$.
We obtain a continuous, even smooth,  finite $\Gamma$-energy $\Gamma$-equivariant mapping 
$$
g: U_X \longrightarrow (H^b)^*\qquad \big(u \mapsto (\phi_0(u), \phi_1(u), \dots)\big).
$$
Also we get a natural mapping $g: U_X \longrightarrow \bold P((H^b)^*),\;u
\mapsto \psi(u)\, (\forall \psi \in H^b)$. 


We assume $g$ is harmonic; otherwise, we replace $g$ by a harmonic mapping
homotopic to g.
Let $\bold F_\bold C(\infty,0)$ denote the complex flat Fubini space, i.e.  a complex Hilbert space. 

\enddemo

\proclaim{Lemma 1} With assumptions of the theorem, $g$ will produce a pluriharmonic mapping $g^{fl}$. There exists a natural holomorphic mapping $g^h : U_X \longrightarrow \bold F_\bold C(\infty,0)$.
\endproclaim

\demo{Proof of Lemma 1}  We define a harmonic  $\Gamma$-equivariant mapping
$$
g^{fl}: =\eurb S_{g(\bold Q)}\cdot  g: U_X \longrightarrow (H^b)^*.
$$
 We have applied the mapping $g$ followed by 
 the Calabi {\it flattening out}\/ $\eurb S_{g(\bold Q)}$ (a generalized stereographic projection from   $g(\bold Q)$) of the real projective space $\bold F_\bold R(\infty,1)$  
into the Hilbert space \cite{2, Chap.\;4, p.\;17}.
By  \cite{2, Chap.\;4, Cor.\;1, p.\;20}, the whole  $\bold F_\bold R(\infty,1)$,
except the antipolar hyperplane $A$ of  $g(\bold Q)$, can be flatten out into $\bold F_\bold R (\infty,0)$. The image of $g$ does not intersect the antipolar hyperplane $A$ of $g(\bold Q)$. Thus we have introduced a flat metric in a large (i.e.\;outside $A$)  neighborhood  of  $g(\bold Q)$ in $ \bold P((H^b)^*)$.

Since the mapping $g^{fl}$ has finite $\Gamma$-energy, it is pluriharmonic; this  is a  special case of a theorem of Siu (see, e.g., \cite{1}). 
Since $U_X$ is simply connected, we obtain the natural {\it holomorphic}\/ mappings 
$$
g^h:
U_X \longrightarrow \bold F_\bold C (\infty,0) 
\big(\hookrightarrow \bold P_\bold C
((H^b)^*)=\bold F_\bold C(\infty,1)\big).
$$
\enddemo

\proclaim{Lemma 2} Construction of a complex line bundle $\Cal L_X$ on $X$ and 
its pullback on $U_X$, denoted by $\Cal L$.
\endproclaim

\demo{Proof of Lemma 2} We take a  point $u\in U_X$. Let $v :=
g^h(u) \in \bold F_\bold C(\infty, 1)$, where $ \bold F_\bold C(\infty, 1)$ is the complex projective space. We consider the linear system of hyperplanes in
$\bold F_\bold C(\infty, 1)$ through $v$ and its proper transform on $U_X$. 
We consider only the moving part. The projection
on $X$ of the latter linear system on $U_X$ will produce a linear system on $X$. 

A connected component of
a {\it general} member of the latter linear system on $X$ will be an irreducible divisor 
$D$ on $X$ by Bertini's theorem. The corresponding line bundle will be the desired
$\Cal L_X:=\Cal O_X(D)$ on $X$. 
\enddemo

\proclaim{Lemma 3} Conclusion of the proof of theorem by induction on $\dim X$.
\endproclaim

\demo{Proof of Lemma 3} 
By the Campana-Deligne theorem \cite{5, Theorem 2.14}, $\pi_1(D)$ will be nonamenable. 
We proceed by induction on $\dim X$, the case
$\dim X=1$ being trivial. Let $q=q(n)$ be an appropriate integer.

We get a  global holomorphic function-section  
$f$ of $\Cal L^q$ corresponding to a bounded pluriharmonic function (see Lemma 1 
and \cite{9, Sect.\;4}).
We will define a $\Gamma$-invariant Hermitian quasi-metric  on sections of $\Cal L^q$  below.
Furthermore,  $f$  is $\ell^2$ on orbits of $\Gamma$, and it  is not  identically zero on any orbit because, otherwise, we could have replaced
$U_X$ by $U_X\backslash B$, where the closed analytic subset $B\subset U_X$ is the union of those orbits on which $f$ had vanished \cite{5, Theorem 13.2, Proof of
Theorem 13.9}.

One  can show that  $f$ satisfies the above conditions by taking linear systems of curvilinear sections of $U_X$ through  $u\in U_X$ and their projections on $X$ (see the proof of Lemma 2 above), since  the statements are trivial in dimension one.
The required Hermitian quasi-metric on $\Cal L^q_X$ is also defined by 
induction on dimension with the help of the Poincar\'e residue map \cite{3,
pp.\;147-148}.

The condition $\ell^2$ on orbits of $\Gamma$ is a local property on $X$. We
get only a Hermitian quasi-metric on $\Cal L^q_X$ (instead of a Hermitian metric). Precisely, 
we get Hermitian metrics over small neighborhoods of  points of $X$, and on the intersections of neighborhoods, they will differ by constant multiples (see \cite{5, Chap.\;5.13}).

For $\forall k>N\gg 0$, the Poincar\'e series are continuous sections
$$
P(f^k)(u):= \sum_{\gamma\in \Gamma} \gamma^*f^k(\gamma u),
$$
and they do not vanish for infinitely many $k$ (see \cite{5, Sect.\;13.1, Theorem 13.2}).

Finally, we can apply Gromov's theorem, precisely, its generalization by Koll\'ar 
(see \cite{4, Corollary 3.2.B, Remark 3.2.B$'$} and 
\cite{5, Theorem 13.8, Corollary 13.8.2, Theorem 13.9, Theorem 13.10}).
So, $X$ is a Moishezon manifold.

The Lemma 3 and Theorem are established.

\enddemo

\remark{Remarks} i) The theorem of the present note provides an alternative proof of a conjecture of H. Wu  provided $\pi_1(X)$ is residually finite (see \cite{10}).

ii) A generalization of the theorem to singular spaces will appear elsewhere.


\endremark

\enddocument